\documentclass[11pt,a4paper]{article}
\usepackage{fontspec}
\usepackage[english]{babel}
\usepackage{amsmath,amssymb,amsthm,mathtools}
\usepackage[a4paper,top=27mm,bottom=27mm,left=30mm,right=30mm]{geometry}
\usepackage{xurl,xcolor}
\usepackage{csquotes}
\usepackage[backend=biber,style=alphabetic,maxbibnames=99,doi=false,url=false]{biblatex}
\definecolor{paperblue}{HTML}{25435B}
\usepackage[unicode,pdftitle={Kähler groups are projective},pdfauthor={Tianshu Jiang and Jiayu Li},colorlinks=true,linkcolor=paperblue,citecolor=paperblue,urlcolor=paperblue]{hyperref}
\DeclareFieldFormat[article,incollection,thesis,unpublished,misc,online]{title}{%
  \mkbibemph{\iffieldundef{doi}
    {\iffieldundef{url}{#1}{\href{\thefield{url}}{#1}}}
    {\href{https://doi.org/\thefield{doi}}{#1}}}}
\renewbibmacro*{in:}{\ifentrytype{article}{}{\printtext{\bibstring{in}\intitlepunct}}}

\AtBeginBibliography{\raggedright}
\theoremstyle{plain}
\newtheorem{theorem}{Theorem}[section]
\newtheorem{lemma}[theorem]{Lemma}
\newtheorem{proposition}[theorem]{Proposition}

\newtheorem*{reltheorem}{Theorem A}
\newtheorem*{projectivecorollary}{Corollary B}
\theoremstyle{definition}

\numberwithin{equation}{section}
\title{Projective Realizations of Kähler Groups}
\author{Tianshu Jiang and Jiayu Li}
\date{}
\begin{document}
\maketitle
\begin{abstract}
We prove that every Kähler group is the fundamental group of a smooth complex projective variety. The main step concerns a surjective holomorphic map \(f\colon M\to B\) with connected fibres from a connected compact Kähler manifold to a smooth projective variety. If the fundamental group of a general smooth fibre has virtually abelian image in \(\pi_1(M)\), then \(\pi_1(M)\) is projective. For the algebraic reduction, a very general smooth fibre is special, and the abelianity theorem for special compact Kähler manifolds gives the required condition on its fundamental group.
\end{abstract}

\section{Introduction}\label{sec:reduction}

A group is called \emph{Kähler} if it is the fundamental group of a connected compact Kähler manifold, and \emph{projective} if it is the fundamental group of a smooth connected complex projective variety. Johnson and Rees asked whether every Kähler group is projective \cite[p.~463]{johnsonrees1987}. Johnson answered the question for compact flat Kähler manifolds \cite[Theorem~I]{johnson1991}, and Baues and Riesterer proved that every virtually abelian Kähler group is projective \cite[Theorem~1.4]{baues2011}.

For a Kähler group admitting a faithful finite-dimensional complex representation, Campana, Claudon and Eyssidieux first obtained a projective subgroup of finite index \cite[Theorem~0.2]{cce2014}. They also showed that a smooth family of complex tori with compact Kähler total space over a projective base can be deformed over that base to a family with projective total space \cite[Theorem~1]{cceadd2015}. Claudon constructed such deformations equivariantly under a finite group action, with projective total spaces for a dense set of parameters. He deduced that every complex linear Kähler group is itself projective \cite[Theorem~1.2 and Corollaries~1.3--1.4]{claudon2018}.

The torus fibrations in Claudon's proof suggest the following relative form of the theorem of Baues and Riesterer.

\begin{reltheorem}
Let \(f:M\to B\) be a surjective holomorphic map with connected fibres from a connected compact Kähler manifold to a smooth connected complex projective variety. For a general smooth fibre \(F\), put
\[
 G=\pi_1(M),\qquad K=\operatorname{im}\bigl(\pi_1(F)\longrightarrow G\bigr).
\]
If \(K\) is virtually abelian, there exists a smooth connected complex projective variety \(Y\) with \(\pi_1(Y)\simeq G\).
\end{reltheorem}

Let \(X\) be a connected compact Kähler manifold. Choose its algebraic reduction \(X\dashrightarrow B\) with \(B\) smooth and projective. Resolving the indeterminacy gives a holomorphic map \(f:M\to B\) with connected fibres from a compact Kähler modification \(M\) of \(X\). By Campana's theorem, a very general smooth fibre \(F\) is special \cite[Theorem~2.39]{campana2004}. The abelianity theorem of OpenAI \cite[Theorem~1.1]{openai2026} then implies that \(\pi_1(F)\) is virtually abelian, and hence so is its image in \(\pi_1(M)\). Theorem~A therefore gives a smooth projective realization of \(\pi_1(M)\simeq\pi_1(X)\).

\begin{projectivecorollary}
For every connected compact Kähler manifold \(X\), there exists a smooth connected complex projective variety \(Y\) with \(\pi_1(Y)\simeq\pi_1(X)\).
\end{projectivecorollary}

Claudon, Höring and Lin proved Corollary~B for compact Kähler threefolds and, in every dimension, for compact Kähler manifolds of algebraic dimension \(a(X)=\dim X-1\) \cite[Theorem~1.2 and Corollary~1.5]{claudon2019}. The algebraic reduction already enters their proof in dimension three \cite[\S4B]{claudon2019}. In algebraic codimension one they found a smooth bimeromorphic model admitting an algebraic approximation \cite[Theorem~1.4]{claudon2019}. In both settings Lin realizes \(\pi_1(X)\) by arbitrarily small projective deformations of \(X\) itself \cite[Theorem~1.1]{linthreefolds}, \cite[Theorem~1.1]{lin2020}. 

To realize \(\pi_1(M)\) projectively in Theorem~A, we start from the fundamental group sequence of the fibration \(f\colon M\to B\). We first identify the quotient of \(G=\pi_1(M)\) by the virtually abelian image \(K\) of the general fibre group. We modify the source and base bimeromorphically as in Proposition~\ref{prop:neat}, preserving \(G\) and its subgroup \(K\). For each divisor in the discriminant of the modified fibration, the orbifold order is the greatest common divisor of the multiplicities of the fibre components over its general point. The resulting orbifold base \(\mathcal B\) has fundamental group \(G/K\) by the exact sequence of the fibration \cite[Prop.~11.9 of the preprint]{camporb2011}. We choose a characteristic free abelian subgroup \(A\subset K\) of finite index and construct a finite gerbe over \(\mathcal B\) whose fundamental group is \(G/A\). The lifted fibration has fibre image \(A\).

Over a dense open substack \(\mathcal U\) the lifted fibration is smooth. In the first real cohomology of a fibre, \(\operatorname{Hom}(A,\mathbb R)\) is the subspace fixed by the kernel of \(\pi_1(\mathcal U)\to G/A\). These subspaces form a variation of Hodge structure. We extend it across the complement of \(\mathcal U\) and pull it back to a smooth projective variety with fundamental group \(G/A\). Deforming its Hodge filtration gives an abelian scheme whose total space has fundamental group \(A\rtimes(G/A)\). A crossed homomorphism \(G\to A\otimes\mathbb R\) restricts to the natural inclusion on \(A\). Choosing such a homomorphism over \(\mathbb Q\) and clearing denominators embeds \(G\) as a subgroup of finite index in \(A\rtimes(G/A)\). The corresponding finite cover of the total space realizes \(G\).

\section{Multiple fibres and fundamental groups of stacks}\label{sec:preparation}

Let \(f\colon X\to B\) be a surjective holomorphic map with connected fibres from a connected compact Kähler manifold to a smooth connected projective variety. Fix a base point on a general smooth fibre \(F\), and write
\[
 G=\pi_1(X),\qquad K=\operatorname{im}\bigl(\pi_1(F)\longrightarrow G\bigr).
\]
We use the classical multiplicities of the orbifold base, given by the greatest common divisors of the fibre multiplicities. Campana proved the corresponding exact sequence for a neat fibration \cite[Prop.~11.9 of the preprint]{camporb2011}. Arapura gave the local lift to the orbifold under normal crossing hypotheses \cite[Lemma~3.5]{arapura2011}. The proof below uses meridians at general points of the components of the inverse image of the discriminant.

\begin{proposition}\label{prop:neat}
There is a bimeromorphic model \(f'\colon M\to B'\) with \(M\) smooth compact Kähler and \(B'\) smooth projective, and a reduced simple normal crossings divisor \(\Delta=\sum_i\Delta_i\) on \(B'\), with the following properties. Let \(m_i\) be the greatest common divisor of the multiplicities of the components of \(f'{}^*\Delta_i\) that dominate \(\Delta_i\). Form the root stack \(\mathcal B\) of the divisors \(\Delta_i\) with these orders. There is a closed analytic subset \(Z\subset B'\) of codimension at least two such that, on writing
\[
 \mathcal B^\circ=\mathcal B|_{B'\setminus Z},\qquad M^\circ=M\setminus f'{}^{-1}(Z),
\]
the map admits a proper representable surjective holomorphic lift \(p\colon M^\circ\to(\mathcal B^\circ)^{\mathrm{an}}\). Under the natural identifications \(\pi_1(M^\circ)\simeq G\) and \(\pi_1(\mathcal B^\circ)\simeq\pi_1(\mathcal B)\), this lift induces an exact sequence
\[
 1\longrightarrow K\longrightarrow G\longrightarrow \pi_1(\mathcal B)\longrightarrow1 .
\]
Over \(U=B'\setminus(Z\cup\Delta)\) the family is smooth with connected fibres. A Kähler form on \(M\) induces compatible Kähler forms on the sources obtained by pulling back to étale manifold charts of the target.
\end{proposition}

\begin{proof}
If \(B\) is a point, take \(M=X\), \(B'=B\), and \(Z=\Delta=\varnothing\); then \(K=G\). Suppose that \(\dim B>0\).

Choose a neat model as in \cite[Definition~1.2 and Lemma~1.3]{campana2004}, with the smooth target of the source modification prescribed to be \(X\). Thus there is a modification \(u\colon M\to X\) such that every prime divisor \(E\subset M\) with \(\operatorname{codim}_{B'}f'(E)\geq2\) has \(\operatorname{codim}_X u(E)\geq2\). In the flattening construction the modifications of the base are blowups; resolving the base and its discriminant by further blowups leaves \(B'\) projective. On the source, take an embedded resolution of the component of the graph dominating \(X\) in \(X\times B'\). It is a smooth submanifold of a succession of blowups of a compact Kähler manifold along smooth centres, and is therefore compact Kähler. Further smooth models dominating a neat model remain neat: a divisor either is contracted by the source modification or maps to an exceptional divisor of the earlier fibration. All these modifications may be chosen to be isomorphisms over a common dense smooth open set. The resulting map has connected fibres, also by Stein factorization over the normal base \(B'\).

Take \(\Delta\) to contain the divisorial discriminant. Let \(Z\) contain the remaining critical-value components, the singular locus of \(\Delta\), and the images of all \(f'\)-exceptional prime divisors. These are finitely many closed analytic subsets of codimension at least two. Since \(B'\) is projective, Chow's theorem makes these subsets algebraic \cite[Theorem~VI]{chow1949}. In particular, no divisor of \(M^\circ\) maps to codimension at least two in \(B'\setminus Z\).

We first check the fundamental group of the source. Deleting a proper closed analytic subset of a complex manifold induces a surjection on fundamental groups; deleting one of complex codimension at least two induces an isomorphism. Perturbing paths, loops and homotopy discs relative to their boundaries to be transverse to a stratification proves these assertions. Let \(C_0\subset X\) be a closed analytic subset of codimension at least two over whose complement \(u\) is an isomorphism. Every divisorial component of \(f'{}^{-1}(Z)\) is contracted by neatness; every remaining component already has dimension at most \(\dim X-2\), as does its image under \(u\). Properness of \(u\) therefore shows that
\[
 C=C_0\cup u(f'{}^{-1}(Z))
\]
is closed analytic of codimension at least two in \(X\). The inclusions
\[
 X\setminus C\simeq u^{-1}(X\setminus C)\subset M^\circ\subset M
\]
induce surjections on fundamental groups. Composing either of the maps from \(X\setminus C\) to \(M^\circ\) or \(M\) with \(u\) gives the inclusion into \(X\), which induces an isomorphism. Hence
\(\pi_1(M^\circ)\simeq\pi_1(M)\simeq G\).
The common smooth fibre gives exactly the original subgroup \(K\) under these identifications; changing the general base point only conjugates it.

Every component of \((f'|_{M^\circ})^*\Delta_i\) dominates \(\Delta_i\setminus Z\), so its multiplicity is divisible by \(m_i\). Consequently
\[
 D_i=\frac{1}{m_i}(f'|_{M^\circ})^*\Delta_i
\]
is an integral Cartier divisor on the smooth manifold \(M^\circ\). The line bundle \(\mathcal O(D_i)\), its canonical section, and the identification of its \(m_i\)-th power with the pulled-back divisor pair define \(p\) by the universal property of the root stack \cite[\S2.2]{cadman2007}. The same description holds holomorphically. The morphism is representable because a manifold has trivial inertia. It factors as
\[
 M^\circ\longrightarrow
 M^\circ\times_{B'\setminus Z}(\mathcal B^\circ)^{\mathrm{an}}
 \longrightarrow (\mathcal B^\circ)^{\mathrm{an}}.
\]
The first arrow is a base change of the finite relative diagonal of the root stack; the second is a base change of the proper map \(f'|_{M^\circ}\). Thus \(p\) is proper. It is surjective: over every geometric point of the root stack, an original nonempty fibre contains a point whose root data are isomorphic to those of that geometric point. For an étale manifold chart \(V\to(\mathcal B^\circ)^{\mathrm{an}}\), the source \(M^\circ\times_{\mathcal B^\circ}V\) is étale over \(M^\circ\), hence smooth. The pullback of one chosen Kähler form \(\omega\) on \(M\) is Kähler there. These forms agree under chart transitions because they all pull back \(\omega|_{M^\circ}\).

Put \(P=\pi_1(U)\), and let \(\gamma_i\) be a positive meridian around \(\Delta_i\). Deleting \(Z\) from \(B'\setminus\Delta\) does not change its fundamental group. In the root chart transverse to \(\Delta_i\), filling the punctured disc by \([\mathbb D/\mu_{m_i}]\) imposes the relation \(\gamma_i^{m_i}=1\); at crossings the charts are products of these root charts. Van Kampen consequently gives compatible presentations
\[
 \pi_1(\mathcal B^\circ)\simeq\pi_1(\mathcal B)
 \simeq P/\langle\!\langle\gamma_i^{m_i}\ (i)\rangle\!\rangle .
\]
This argument computes the fundamental group of the stack classifying space and does not assume that inertia injects into it.

For completeness, let \(H=\pi_1(f'{}^{-1}U)\). The homotopy sequence of the smooth proper family over \(U\) gives a surjection \(\psi\colon H\to P\), whose kernel \(N_F\) is the image of \(\pi_1(F)\) in \(H\). The kernel \(R\) of \(H\to\pi_1(M^\circ)=G\) is normally generated by the meridians of the vertical divisors. Indeed, a bounding disc can be made transverse to their smooth parts and disjoint from their singular and intersection loci; removing small discs around its intersection points gives this meridian description. If a component over \(\Delta_i\) has multiplicity \(n_{ij}\), its meridian maps under \(\psi\) to a conjugate of \(\gamma_i^{n_{ij}}\). Since \(\psi\) is surjective, \(\psi(R)\) is the normal closure of all these powers in \(P\). Divisibility and a Bézout identity for the integers \(n_{ij}\) show that
\[
 \psi(R)=\langle\!\langle\gamma_i^{m_i}\ (i)\rangle\!\rangle.
\]
The image of \(N_F\) in \(H/R=G\) is precisely \(K\). Therefore
\[
 G/K=H/(N_FR)=P/\psi(R)=\pi_1(\mathcal B).
\]
This quotient is induced by \(p\), since the maps agree over \(U\), proving the asserted exact sequence.
\end{proof}

Hu and Wu proved the following assertion for effective stacks \cite[Lemma~6]{huwu2026}. Their proof uses the construction in \cite[Lemma~7.1]{openai2026b}. The gerbe in Section~\ref{sec:transfer} can have nontrivial generic stabilizer, so we also need the statement for such stacks.

\begin{proposition}\label{prop:projective-model}
Let \(\mathcal T\) be a smooth connected proper complex algebraic Deligne--Mumford stack with projective coarse space \(T\). There is a smooth connected projective variety \(S\) and an algebraic morphism \(S\to\mathcal T\) inducing an isomorphism on fundamental groups.
\end{proposition}

\begin{proof}
The stack \(\mathcal T\) is separated and tame and has projective coarse space. By the quotient criterion of Kresch and Vistoli \cite[Theorem~2]{kv2004} and Gabber's theorem on the Brauer group \cite[Theorem~1.1]{dejong}, \(\mathcal T\) is a quotient by a linear algebraic group. In dimension zero this follows directly from \(\mathcal T\simeq\mathrm BH\) for a finite group \(H\). By \cite[Lemma~2.12]{ehkv2001}, there is a vector bundle \(V\) on which every geometric stabilizer acts faithfully. Adding a trivial summand if necessary, take \(r=\operatorname{rk}V\geq1\). This includes the ineffective stabilizers of \(\mathcal T\).

Put \(\mathcal P=\mathbb P_{\mathcal T}(\mathcal O\oplus V)\), with projectivization parametrizing lines. Every nonidentity stabilizer element acts nonscalarly on \(\mathcal O\oplus V\): it fixes the trivial line and acts nontrivially on \(V\). Its projective fixed locus is a union of proper projective eigenspaces. Since the stabilizer is finite, each projective fibre contains a point with trivial stabilizer. The free locus of \(\mathcal P\) is open: the complement of the identity section in its finite inertia has closed image. The smooth connected stack \(\mathcal P\) is irreducible, so this nonempty open locus is dense. A finite group acting on a connected smooth local chart cannot act identically on a nonempty open subset without acting identically on the chart. Thus trivial generic stabilizer implies effectiveness here.

To see that the coarse space of \(\mathcal P\) is projective, we use the descent argument of \cite[Lemma~2 and proof of Theorem~1]{kv2004}. Finite inertia on a finite-type stack bounds the stabilizer orders. For a common multiple \(\ell\) of these orders, every stabilizer acts trivially on the fibres of \(\mathcal O_{\mathcal P}(\ell)\), so this line bundle descends to the coarse space. Descent can be checked on finite quotient charts: a linearized line bundle whose point stabilizers act trivially has invariant local generators, obtained by averaging, and these generators descend. On an étale chart of the coarse base, write \(\mathcal T=[A/\Gamma]\) with \(\Gamma\) finite and \(A\) finite over that chart. The coarse projective bundle is the finite quotient of \(\mathbb P_A((\mathcal O\oplus V)|_A)\). The pullback of the descended line bundle is relatively ample; ampleness descends along a finite surjective map. Hence the descended line bundle is relatively ample over \(T\). Twisting a positive power by a sufficiently ample pullback from \(T\) proves that the coarse space is projective.

The projection \(\mathcal P\to\mathcal T\) induces an isomorphism on fundamental groups. On classifying spaces it is a projective-space bundle, whose fibre \(\mathbb P^r(\mathbb C)\) is connected and simply connected. By \cite[Lemma~6]{huwu2026}, there is a morphism from a smooth connected projective variety \(S\) to \(\mathcal P\) inducing an isomorphism on fundamental groups. Composing with the projection to \(\mathcal T\) gives the required morphism.
\end{proof}

\section{Finite gerbes}\label{sec:transfer}

Suppose that the image \(K\) of the fundamental group of a general smooth fibre \(F\) is virtually abelian. Choose a characteristic free abelian subgroup \(A\subset K\) of finite index, and put \(N=K/A\). We construct a gerbe over the base that is étale locally a product with \(\mathrm BN\). The original map lifts to this gerbe, and a general geometric fibre of the lift is a connected finite cover of \(F\).

\begin{proposition}\label{prop:transfer}
Let \(\mathcal B\) be a smooth connected proper complex algebraic Deligne--Mumford stack, and let \(\mathcal B^\circ\hookrightarrow\mathcal B\) be a connected open substack inducing an isomorphism on fundamental groups. Let \(p\colon M\to(\mathcal B^\circ)^{\mathrm{an}}\) be a proper representable surjective holomorphic map from a smooth connected Kähler manifold. Suppose that \(p\) is smooth with connected geometric fibres over a nonempty connected dense open substack, and that it induces an exact sequence
\[
 1\longrightarrow K\longrightarrow G=\pi_1(M)
   \xrightarrow{p_*} Q=\pi_1(\mathcal B)\longrightarrow1 ,
\]
where \(K\) is the image in \(G\) of the fundamental group of a general geometric fibre \(F\).
Let \(A\trianglelefteq G\) satisfy \(A\subset K\) and \([K:A]<\infty\), and set \(N=K/A\) and \(E=G/A\). There exists a proper algebraic gerbe \(a\colon\mathcal E\to\mathcal B\) that is étale locally a product with \(\mathrm BN\). The stack \(\mathcal E\) is smooth and connected and has the same coarse space as \(\mathcal B\). For \(\mathcal E^\circ=\mathcal E\times_{\mathcal B}\mathcal B^\circ\) we have
\[
 \pi_1(\mathcal E^\circ)\simeq\pi_1(\mathcal E)\simeq E,\qquad
 a_*=(E\twoheadrightarrow Q).
\]
There is a proper representable surjective holomorphic lift \(p'\colon M\to(\mathcal E^\circ)^{\mathrm{an}}\) inducing \(G\twoheadrightarrow E\). It is smooth over the preimage of the original smooth open substack. Its general geometric fibre is the connected cover of \(F\) classified by \(\pi_1(F)\twoheadrightarrow N\), and its fundamental group has image \(A\) in \(G\).
\end{proposition}

\begin{proof}
Write \(\epsilon\colon E\twoheadrightarrow Q\) for the quotient map. For a discrete group \(H\), let \(\mathrm BH\) denote the classifying stack of principal \(H\)-coverings. A smooth analytic Deligne--Mumford stack has local charts by balls with finite group actions. Hence \(\mathcal B^{\mathrm{an}}\) is locally path connected and semilocally simply connected in the sense of \cite[Theorem~18.19 and Corollary~18.20]{noohi2005}. It has a universal principal \(Q\)-cover with simply connected total stack. A covering inherits a unique compatible complex structure on each manifold chart, so its classifying map \(c\colon\mathcal B^{\mathrm{an}}\to\mathrm BQ\) is holomorphic. Form
\[
 \mathcal G=\mathcal B^{\mathrm{an}}\times_{\mathrm BQ}\mathrm BE .
\]
Over a chart trivializing the \(Q\)-cover, \(\mathcal G\) is a product with \(\mathrm BN\). The principal \(E\)-cover pulled back from \(\mathrm BE\) has total stack \(\mathcal B^{\mathrm{an}}\times_{\mathrm BQ}*\), the same simply connected universal cover. Thus \(\pi_1(\mathcal G)=E\), with projection \(\epsilon\) to \(Q\). Relative inertia acts by multiplication on the fibres of this \(E\)-cover, giving the specified inclusion \(N\hookrightarrow E\). The restriction of the universal cover to \(\mathcal B^\circ\) is again connected and simply connected because \(\pi_1(\mathcal B^\circ)\to Q\) is an isomorphism.

To algebraize \(\mathcal G\), let \(\rho\colon Q\to\operatorname{Out}(N)\) be the outer action induced by conjugation in \(E\). Write \(Z_N=Z(N)\) and factor \(\epsilon\) as
\[
 E\xrightarrow{\ \beta\ } E_0=Q\times_{\operatorname{Out}(N)}\operatorname{Aut}(N)
   \longrightarrow Q,\qquad
 \beta(e)=(\epsilon(e),\operatorname{Ad}(e)|_N).
\]
Two representatives of the same outer automorphism differ by an inner automorphism. Thus \(\beta\) is surjective with kernel \(Z_N\). The finite principal \(\operatorname{Out}(N)\)-cover defined by \(\rho\) is algebraic by the Riemann existence theorem for stacks \cite[Theorem~20.4]{noohi2005}. Its classifying morphism defines an algebraic stack
\[
 \mathcal C=\mathcal B\times_{\mathrm B\operatorname{Out}(N)}
                       \mathrm B\operatorname{Aut}(N).
\]
This is a gerbe locally modelled on \(\mathrm B(N/Z_N)\). Its analytification is \(\mathcal B^{\mathrm{an}}\times_{\mathrm BQ}\mathrm BE_0\). Indeed, taking the fibre product of the \(Q\)-cover and the lifted \(\operatorname{Aut}(N)\)-cover over their common \(\operatorname{Out}(N)\)-cover gives the required \(E_0\)-cover. The map \(\mathcal G\to\mathcal C^{\mathrm{an}}\) is therefore an abelian gerbe banded by \(\mathcal Z^{\mathrm{an}}\). Here \(\mathcal Z\) is the finite algebraic local system obtained from the natural action of \(\operatorname{Aut}(N)\) on \(Z_N\). This allows nontrivial monodromy on the centre.

The required comparison is
\[
 H^2_{\text{ét}}(\mathcal C,\mathcal Z)
       \simeq H^2(\mathcal C^{\mathrm{an}},\mathcal Z^{\mathrm{an}}).
\]
Choose a surjective étale atlas that is a finite disjoint union of smooth affine varieties. The diagonal of \(\mathcal C\) is finite, so every term of the Čech nerve is again smooth and affine. Artin's comparison theorem \cite[Exposé~XI, Théorème~4.4(iii)]{artinsga} applies in each degree of this nerve. The two descent spectral sequences give the displayed isomorphism. By Giraud's classification of gerbes with a fixed abelian band \cite[Ch.~III, Def.~3.1.1 and Thm.~3.3.3]{giraud1966}, the class of \(\mathcal G\to\mathcal C^{\mathrm{an}}\) defines an algebraic gerbe \(\mathcal E\to\mathcal C\) whose analytification is equivalent to \(\mathcal G\) over \(\mathcal C^{\mathrm{an}}\). This gerbe is an algebraic stack because étale locally it is the classifying stack of the finite étale band. For the composite \(\mathcal E\to\mathcal B\), objects and isomorphisms exist étale locally by lifting successively through the two gerbes. The relative automorphism group of a local object is finite étale, and the analytic equivalence identifies its geometric fibres with \(N\). Trivializing that group gives the local product with \(\mathrm BN\). Smoothness, properness, connectedness and the assertion about the coarse space follow from this description.

The quotient \(G\twoheadrightarrow E\) defines a holomorphic principal \(E\)-cover on \(M\). Its associated \(Q\)-cover and the pullback of the universal cover along \(p\) have the same based monodromy. An isomorphism between them gives an object of \(p^*\mathcal G\), hence the lift \(p'\). Its classifying map induces the prescribed quotient \(G\to E\). The lift factors through the section of the pulled-back gerbe on \(M\), followed by a base change of \(p\). Pulling the section back to a manifold gives a finite étale sheaf of isomorphisms. The lift is therefore representable and proper, and is smooth where \(p\) was smooth. It is also surjective: over a geometric object of the gerbe, its fibre is a principal finite \(N\)-cover of the corresponding nonempty fibre of \(p\). Pullback of the Kähler form of \(M\) again gives compatible Kähler forms on the sources of its étale charts.

Fix an object of the gerbe at a general point. The new geometric fibre is a principal \(N\)-cover \(\widehat F\to F\). Its monodromy, composed with the inertia inclusion \(N\hookrightarrow E\), is the composite \(\pi_1(F)\to G\to E\). Thus the monodromy is \(\pi_1(F)\twoheadrightarrow K\twoheadrightarrow N\) and \(\widehat F\) is connected. A loop on \(F\) lifts to a loop on \(\widehat F\) exactly when its monodromy in \(N\) is trivial. Since \(\pi_1(F)\to K\) is surjective, the images of these loops in \(G\) form precisely \(A\).
\end{proof}

If \(K\) is finitely generated and virtually abelian, choose a finite-index abelian subgroup and multiply it by the exponent of its torsion subgroup. The result is a finite-index free abelian subgroup \(K_1\) of \(K\). Put \(d=[K:K_1]\). There are only finitely many subgroups of index \(d\) in a finitely generated group. Their intersection \(A\) is characteristic of finite index in \(K\) and is contained in \(K_1\). Hence \(A\) is free abelian and normal in \(G\). The lift constructed in Proposition~\ref{prop:transfer} induces the exact sequence
\[
 1\longrightarrow A\longrightarrow G\longrightarrow\pi_1(\mathcal E)=E\longrightarrow1 .
\]
\section{Fibre cohomology and group extensions}\label{sec:cohomology}

Let \(p\colon M\to\mathcal Y\) be a proper representable surjective holomorphic map from a connected complex manifold to a smooth connected separated complex analytic Deligne--Mumford stack with finite stabilizers. Fix a K\"ahler form \(\omega\) on \(M\). Suppose that \(p\) is smooth with connected geometric fibres over a dense connected Zariski-open substack \(j\colon\mathcal U\hookrightarrow\mathcal Y\); here Zariski-open means that \(\mathcal Y\setminus\mathcal U\) is a closed complex-analytic substack of positive codimension. Choose a geometric fibre \(F\) over \(u\in\mathcal U\) and compatible base points. Fundamental groups of stacks refer to their analytic classifying spaces. Write
\[
 P=\pi_1(\mathcal U),\qquad Q=\pi_1(\mathcal Y),\qquad G=\pi_1(M),\qquad
 K=\operatorname{im}\bigl(\pi_1(F)\to G\bigr).
\]
Assume that \(P\to Q\) is surjective and that \(p_*\) induces an exact sequence \(1\to K\to G\to Q\to1\). Conjugation gives a \(Q\)-action on \(K^{\mathrm{ab}}\). The real dual of its torsion-free quotient is \(\operatorname{Hom}(K,\mathbb R)\). We identify this space inside \(H^1(F,\mathbb R)\) and lift the corresponding evaluation map to a cohomology class on \(M\).

\begin{proposition}\label{prop:cocycle}
Let \(L=\ker(P\to Q)\), \(A=K^{\mathrm{ab}}/\operatorname{Tor}(K^{\mathrm{ab}})\), \(A_{\mathbb R}=A\otimes_{\mathbb Z}\mathbb R\), and \(\mathcal W=R^1(p|_{p^{-1}\mathcal U})_*\mathbb R\).
Pullback along \(\pi_1(F)\twoheadrightarrow K\) gives an isomorphism of \(Q\)-modules
\[
 \operatorname{Hom}(K,\mathbb R)\xrightarrow{\ \sim\ } H^1(F,\mathbb R)^L.
\]
Let \(\mathcal D\) be the local system on \(\mathcal Y\) corresponding to the invariant subspace. There is a class \(c\in H^1(M,p^*\mathcal D^\vee)\) whose restriction to \(F\) corresponds to the inclusion \(\mathcal D_u\hookrightarrow H^1(F,\mathbb R)\). Under the resulting identification \(\mathcal D_u^\vee=A_{\mathbb R}\), a crossed homomorphism
\[
 t\colon G\to A_{\mathbb R},\qquad t(gh)=t(g)+p_*(g)t(h)
\]
represents \(c\), and \(t|_K\) is the natural map \(K\to A\hookrightarrow A_{\mathbb R}\).
\end{proposition}

\begin{proof}
Set \(C=Rp_*\mathbb R_M\), \(n=\dim M\), and \(b=\dim\mathcal Y\). We first specify the analytic sheaf category and the descent used below. Choose an analytic \'etale manifold atlas \(Z_0\to\mathcal Y\), with \v{C}ech nerve \(Z_\bullet\). Take the homotopy limit of the stable enhanced derived categories of real sheaves on \(Z_r\), with transition functors given by inverse image. Its objects are complexes with homotopy-coherent descent data. Derived sheaf descent identifies this category independently of the atlas; a common refinement gives the same category. We use the full subcategory whose restriction to \(Z_0\) has finite-dimensional cohomology sheaves, constructible for a locally finite complex-analytic stratification, and vanishes outside a finite interval of degrees.

Every structure map of \(Z_\bullet\) is a local biholomorphism. Its inverse image preserves the support and cosupport conditions defining middle perversity and is therefore perverse \(t\)-exact. Perverse truncation commutes canonically with these inverse images: the comparison maps are the unique ones compatible with the truncation adjunctions, and consequently satisfy all the coherence conditions. Applying truncation to a descent datum thus gives a descent datum and its truncation triangle. If one object is in perverse degrees at most zero and another in degrees at least one, the mapping space between them is contractible on every term of the nerve; its limit is contractible as well. These observations prove the truncation and orthogonality axioms on the stack. The resulting \(t\)-structure is bounded on the chosen subcategory, since the charts have the fixed finite dimension \(b\). Pullback to the atlas is \(t\)-exact and detects isomorphisms. Its heart and the intermediate-extension functor have the corresponding \'etale descent.

For an \'etale manifold chart \(Z\to\mathcal Y\), the base change \(p_Z\colon M_Z=M\times_{\mathcal Y}Z\to Z\) is proper. The map \(M_Z\to M\) is a local biholomorphism, so \(M_Z\) is smooth and carries the pullback of \(\omega\). This form is closed and remains Kähler after addition of the pullback of a local Kähler form on \(Z\). Thus \(p_Z\) is a Kähler morphism \cite[Remark~1.2]{saito2022}. Saito's decomposition theorem with constant real coefficients, stated in \cite[Theorem~0.5]{saito1990} and proved in \cite[Theorem~1]{saito2022}, applies to each connected component. It gives relative hard Lefschetz for \(Rp_{Z*}\mathbb R_{M_Z}[n]\) and the decomposition by strict support of its perverse cohomology. These direct images are bounded and analytically constructible. Proper base change identifies them with the restrictions of \(C[n]\), compatibly on the nerve.

Cup product with the single global class \([\omega]\) gives a global morphism \(\ell\colon C[n]\to C[n+2]\). The chartwise Lefschetz operators are its pullbacks. Since pullback detects isomorphisms, we obtain
\[
 \ell^i\colon{}^pH^{-i}(C[n])\xrightarrow{\ \sim\ }{}^pH^i(C[n])
 \qquad(i>0)
\]
on \(\mathcal Y\). Deligne's criterion applies to a bounded object in any triangulated category with a \(t\)-structure, as stated and proved in \cite[(0.3) and \S1]{deligne1994}; see also \cite[Proposition~1.2 and Theorem~1.5]{deligne1968}. Apply it in this global category to \(C[n]\) and shift back. It gives an isomorphism
\[
 \bigoplus_k{}^pH^k(C)[-k]\xrightarrow{\ \sim\ }C
\]
inducing the identity on each perverse cohomology object. This normalization follows by lifting the identity of each cohomology object, as in Deligne's proof on p.~117.

On \(\mathcal U\) the ordinary cohomology sheaves of \(C\) are local systems. A local system shifted by \(b\) is perverse, so
\[
 ({}^pH^{b+1}(C))|_{\mathcal U}=\mathcal W[b].
\]
Set \(I=j_{!*}(\mathcal W[b])[-b]\). The analytic boundary condition allows these intermediate extensions in the constructible category. On each chart Saito's decomposition by strict support splits \({}^pH^{b+1}(C)\) into \(I[b]\) and terms supported on the boundary. No proper-support term meets the smooth open, where the restriction is a shifted local system. An intermediate extension has no perverse subobject or quotient supported on the boundary. Hence there are no morphisms in either direction between \(I[b]\) and a boundary term. An endomorphism of \(I[b]\) that vanishes on \(\mathcal U\) also vanishes. It follows that the inclusion of this full-support summand restricting to the identity of \(\mathcal W[b]\) is unique. These inclusions therefore agree on overlaps and descend in the perverse heart. Shifting by \(-(b+1)\) and composing with the global splitting gives
\[
 \iota\colon I[-1]\longrightarrow C,\qquad
 H^1(\iota|_{\mathcal U})=\operatorname{id}_{\mathcal W}.
\]

Normality of \(L\) makes \(\mathcal W_u^L\) a \(P\)-subrepresentation on which \(L\) acts trivially. Since \(P/L=Q\), it defines a local system \(\mathcal D\) on all of \(\mathcal Y\), with an inclusion \(\mathcal D|_{\mathcal U}\hookrightarrow\mathcal W\). On each smooth chart a local system shifted by the dimension is the intermediate extension of its restriction across a closed analytic subset of positive codimension. Thus
\[
 j_{!*}((\mathcal D|_{\mathcal U})[b])=\mathcal D[b].
\]
Apply the intermediate-extension functor to the inclusion, and then compose with \(\iota\). This gives a morphism \(\mathcal D[-1]\to Rp_*\mathbb R_M\). Derived adjunction and the finite-rank local constancy of \(\mathcal D\) identify its group of morphisms with
\[
 \operatorname{Hom}(\mathcal D[-1],Rp_*\mathbb R_M)
 =\operatorname{Hom}(p^*\mathcal D[-1],\mathbb R_M)
 =H^1(M,p^*\mathcal D^\vee).
\]
Denote the resulting class by \(c\). Proper base change and the normalization of \(\iota\) identify \(c|_F\), under
\[
 H^1(F,\mathcal D_u^\vee)
 =\operatorname{Hom}(\mathcal D_u,H^1(F,\mathbb R)),
\]
with the given inclusion. Only constant real coefficients were used in Saito's theorem; this construction does not require \(\mathcal D\) to be a Hodge subvariation.

Pullback along \(\pi_1(F)\twoheadrightarrow K\) injects \(\operatorname{Hom}(K,\mathbb R)\) into \(H^1(F,\mathbb R)\). The homotopy sequence over \(\mathcal U\) permits a loop of the base to be lifted to \(\pi_1(p^{-1}\mathcal U)\). Transport of fibre loops, after passage to \(G\), agrees with conjugation by this lift. Two lifts differ by a fibre element, whose image lies in \(K\), and inner automorphisms of \(K\) act trivially on \(K^{\mathrm{ab}}\). The injection is therefore \(P\)-equivariant for the action on its source pulled back from \(Q\). Its image lies in \(\mathcal D_u\).

The universal cover of \(M\) is simply connected, so it has zero first cohomology with constant vector-space coefficients. The degree-one edge map for this cover identifies \(H^1(M,p^*\mathcal D^\vee)\) with group cohomology. Represent \(c\) by a crossed homomorphism \(t\colon G\to\mathcal D_u^\vee\). Since \(K\) acts trivially on the coefficients, \(t|_K\) is a homomorphism, and every coboundary restricts to zero there. If \(\gamma\in\pi_1(F)\) has image \(\bar\gamma\in K\), the restriction of \(c\) to \(F\) gives
\[
 \langle t(\bar\gamma),d\rangle=d(\gamma)
 \qquad(d\in\mathcal D_u).
\]
For each \(d\), the map \(k\mapsto\langle t(k),d\rangle\) is a homomorphism \(K\to\mathbb R\) pulling back to \(d\). This proves the reverse inclusion and hence the asserted real isomorphism, with its \(Q\)-action. Under \(\mathcal D_u=\operatorname{Hom}(K,\mathbb R)=A_{\mathbb R}^\vee\), the displayed identity says that \(t(k)\) evaluates each homomorphism on \(k\). Thus \(t|_K\) is precisely \(K\to A\hookrightarrow A_{\mathbb R}\).
\end{proof}

\section{Variations of Hodge structure and abelian schemes}\label{sec:hodge}

Over the smooth locus of \(p\), the first cohomology of the compact K\"ahler fibres carries an integral variation of Hodge structure of types \((1,0)\) and \((0,1)\). If the fibre dimension is \(d\geq1\), cup product with \([\omega]^{d-1}\) followed by fibre integration gives a flat real polarization on first cohomology. For zero-dimensional fibres this cohomology vanishes. We show that the invariant subspace of Proposition~\ref{prop:cocycle} inherits the Hodge filtration and use it to construct an abelian scheme.

\begin{lemma}\label{lem:normal-hodge}
Let \(\mathcal U\) be a smooth connected separated algebraic Deligne--Mumford stack of finite type over \(\mathbb C\), with a dense open substack that is a finite \'etale gerbe over a smooth quasi-projective variety. Let \(\mathcal W\) be a polarizable real variation of pure Hodge structure on \(\mathcal U^{\mathrm{an}}\), choose a geometric base point \(u\), and put \(P=\pi_1(\mathcal U,u)\). For every normal subgroup \(L\triangleleft P\), the local subsystem with fibre \(\mathcal W_u^L\) is a polarizable real subvariation.
\end{lemma}

\begin{proof}
Restrict the gerbe to a dense affine open subset of its base. Choose a connected affine \'etale chart that dominates this open subset. Its composite with the projection to the base is \'etale and quasi-finite of finite type. Embed the chart into the normalization of the base in the function field of the chart and remove the image of the complement. After this shrinking the chart is finite \'etale and surjective over the base. It is also finite \'etale over the gerbe: after pullback to a chart neutralizing the gerbe, the map is an isomorphism torsor finite \'etale over the fibre product of the two charts over the base. Take the cover parametrizing orderings of the geometric fibres of this finite \'etale cover of the gerbe, and then a connected component \(Z\). Projection to the first entry is finite \'etale over the original chart, so \(Z\) is a smooth quasi-projective scheme. The subgroup \(\Gamma\) of the permutation group preserving this component gives a presentation \(\mathcal U_0=[Z/\Gamma]\) of the resulting dense open substack. The finite group action may have a kernel; all descent data, including the action on the variation, are retained.

The variety \(Z\) admits a smooth projective compactification. On \(Z\), the underlying real local system is semisimple, and the global flat sections form a pure Hodge substructure in every fibre. Schmid proved the theorem of the fixed part in \cite[Theorem~7.22]{schmid1973}. Deligne explained its extension to polarizable complex variations, attributed semisimplicity in this generality to Nori, and described the Hodge structure on the algebra of flat endomorphisms \cite[\S\S1.11--1.13]{deligne1987}. Applying these results to the complexification of the real variation and its tensor constructions, then taking the real form, gives the assertions on \(Z\).

For a local subsystem on \(\mathcal U_0\), choose a flat projection onto its pullback to \(Z\) and average the projection over \(\Gamma\). The average takes values in the chosen subsystem and acts as the identity on it, so it descends to a flat projection. The action of \(\Gamma\) on global flat sections respects their Hodge structure, so the invariant sections also form a pure Hodge structure. This proves semisimplicity and the fixed-part statement on \(\mathcal U_0\).

For a local system defined on all of \(\mathcal U\), restriction to \(\mathcal U_0\) preserves global flat sections. Indeed, the complement of a proper closed analytic subset in a polydisc chart is connected, so a locally constant section on that complement extends uniquely in a flat trivialization. This also applies to endomorphism and tensor local systems. In particular, flat projections extend, and evaluation of the extended flat sections remains a morphism of Hodge structures at boundary points, since membership in a Hodge subbundle and its conjugate is a closed condition. Uniqueness gives compatibility with the descent maps.

For an irreducible real \(P\)-module \(V\), normality of \(L\) makes \(V^L\) a \(P\)-submodule, so \(V^L\) is either zero or all of \(V\). This alternative is the same for isomorphic irreducible modules. Thus \(\mathcal W_u^L\) is a sum of whole isotypic components. The projection onto their sum is a central idempotent in the semisimple real algebra
\[
 R=\operatorname{End}_P(\mathcal W_u)
   =H^0\bigl(\mathcal U,\operatorname{End}(\mathcal W)\bigr).
\]
By the fixed-part theorem, \(R\) carries a Hodge structure of weight zero. Composition is a Hodge morphism, so the Hodge torus acts by real algebra automorphisms on \(R\). Its real points \(\mathbb C^\times\) form a connected group and hence fix each of the finitely many primitive central idempotents. Our projection therefore has type \((0,0)\). Its image inherits the Hodge filtration, and the restriction of the polarization gives a positive-definite Hodge metric.
\end{proof}

\begin{proposition}\label{prop:jacobian}
Let \(\mathcal T\) be a smooth connected proper complex algebraic Deligne--Mumford stack with projective coarse space. Let \(\mathcal U\subset\mathcal T\) be a dense connected Zariski-open substack satisfying the gerbe hypothesis of Lemma~\ref{lem:normal-hodge}. Suppose \(P=\pi_1(\mathcal U)\to Q=\pi_1(\mathcal T)\) is surjective and put \(L=\ker(P\to Q)\). Let \(\mathcal W\) be a real-polarizable variation of weight one on \(\mathcal U\), of types \((1,0)\) and \((0,1)\), with a flat integral lattice. Given a \(Q\)-stable full lattice \(\Lambda\) in \(D=\mathcal W_u^L\), the corresponding integral local system, the Hodge filtration on \(D_{\mathbb C}\), and the real polarization extend to \(\mathcal T\).

Let \(h\colon S\to\mathcal T\) be a morphism from a smooth connected projective variety inducing an isomorphism on fundamental groups. There is an abelian scheme \(\pi\colon J\to S\) with smooth projective total space and first integral homology local system \(h^*\Lambda^\vee\). Its zero section identifies
\[
 \pi_1(J)=\Lambda^\vee\rtimes Q.
\]
\end{proposition}

\begin{proof}
By Lemma~\ref{lem:normal-hodge}, \(D\) underlies a real-polarizable subvariation of weight one. Its monodromy factors through \(Q\), so the specified lattice \(\Lambda\) defines an integral local system on all of \(\mathcal T\). We use this lattice, without replacing it by the intersection of \(D\) with the original lattice of \(\mathcal W\). The real polarization is \(Q\)-invariant because \(P\to Q\) is surjective, and therefore extends as a flat form. If the lattice has rank zero, take \(J=S\). Assume henceforth that its rank is \(2g>0\).

Trivialize the integral local system on an \'etale polydisc chart \(B\to\mathcal T\). On the inverse image \(B\setminus Z_B\) of \(\mathcal U\), the flat real polarization is a fixed symplectic form. In a real symplectic basis the period map is a single-valued holomorphic map to the bounded realization of Siegel space,
\[
 B\setminus Z_B\longrightarrow
 \Omega_g=\{Z\in\operatorname{Sym}_g(\mathbb C)\mid I-Z^*Z>0\}.
\]
Its coordinates are bounded and hence extend holomorphically across \(Z_B\). If the extension met the boundary of the convex domain \(\Omega_g\), choose a supporting real hyperplane \(\operatorname{Re}\ell=c\) there, with \(\ell\) complex linear. The real part of \(\ell\) composed with the extension would attain its maximum at an interior point of \(B\), although it is strictly smaller than \(c\) on \(B\setminus Z_B\). The maximum principle excludes this. The extended map thus defines a polarized Hodge filtration, compatible by uniqueness with chart transitions and stabilizer actions. Griffiths transversality is automatic for these two Hodge types.

Pull back to \(S\) and dualize. On \(\mathcal H_{\mathbb Z}=h^*\Lambda^\vee\) this gives a real-polarizable variation of weight \(-1\) and types \((-1,0),(0,-1)\). Campana, Claudon and Eyssidieux established the relative deformation criterion for geometric families \cite[Theorem~1.1 and \S2.3]{cce2014}. Claudon's formulation for variations \cite[Definition~3.1 and Proposition~3.6]{claudon2018} applies here: \(S\) is compact K\"ahler, the variation has weight \(-1\) and a real polarization, and a small deformation is defined on the pullback of the same integral local system to \(S\) times a polydisc. We may thus choose a new filtration on \(\mathcal H_{\mathbb Z}\) with a rational polarization. Multiplying this form by a positive integer clears its denominators in one fibre and hence, by flatness, in every fibre. This gives an integral alternating Riemann form \(b\), without changing the lattice. Our sign convention is \(b(u,Iu)>0\) for nonzero real \(u\), where \(I\) is the induced fibre complex structure. Write \(F^\bullet\) for the new filtration and set
\[
 \mathcal H=\mathcal H_{\mathbb Z}\otimes_{\mathbb Z}\mathcal O_S,\qquad
 J=(\mathcal H/F^0\mathcal H)/\mathcal H_{\mathbb Z}.
\]
The Hodge decomposition identifies \(\mathcal H_{\mathbb R}\) with \(\mathcal H/F^0\mathcal H\) as real vector bundles. Thus \(J\to S\) is a smooth proper holomorphic family of complex tori with a zero section, and its first integral homology local system is the translation lattice \(\mathcal H_{\mathbb Z}\).

To prove projectivity, we construct an integral K\"ahler class. Put \(H=\mathcal H_{\mathbb Z,s}\) and let \(\widetilde S\) be the universal cover of \(S\). As a real manifold, \(J\) is the quotient of \(\widetilde S\times H_{\mathbb R}\) by lattice translations and the local-system action over deck transformations. Let \(v\) denote the fibre coordinate. The constant vertical two-form \(\beta\) with value \(b\), extended by zero on flat horizontal vectors, is closed and invariant under both actions. Each constant real section is a flat holomorphic section of \(\mathcal H\), so its image in \(\mathcal H/F^0\mathcal H\) is holomorphic. The flat horizontal leaves are therefore complex submanifolds. Both the vertical and the horizontal real tangent distributions are invariant under the complex structure of \(J\). The Riemann bilinear relations on the vertical summand then make \(\beta\) a semipositive form of type \((1,1)\), positive in all nonzero vertical directions and with kernel exactly the horizontal distribution.

On the trivial smooth Hermitian line bundle over \(\widetilde S\times H_{\mathbb R}\), let \(a\in H\) act by
\[
 (s,v,z)\longmapsto (s,v+a,\exp(2\pi i\,b(a,v))z).
\]
Deck transformations act on \(v\) by the local-system action and fix \(z\). For \(a,c\in H\), the extra multiplier in composing the two translations is \(\exp(2\pi i\,b(a,c))=1\). Thus these translations form an action, and invariance of \(b\) makes it compatible with the deck action. The quotient is a smooth Hermitian line bundle \(\mathcal L\) on \(J\).

Set \(\theta=b(v,dv)\). Under translation by \(a\), this one-form changes by \(b(a,dv)\), while the logarithmic differential of the multiplier is \(2\pi i\,b(a,dv)\). The unitary connection \(\nabla=d-2\pi i\theta\) consequently descends to \(\mathcal L\). For a flat basis \(e_i\), put \(b_{ij}=b(e_i,e_j)\) and write \(v=\sum_i v_i e_i\). Then \(\beta=\sum_{i<j}b_{ij}\,dv_i\wedge dv_j\) and \(d\theta=2\beta\), so
\[
 F_\nabla=-4\pi i\,\beta,\qquad
 c_1(\mathcal L,\nabla)=\frac{i}{2\pi}F_\nabla=2\beta.
\]
Because this curvature has type \((1,1)\), the operator \(\nabla^{0,1}\) is integrable and defines a holomorphic structure on \(\mathcal L\); \(\nabla\) is its Chern connection. This supplies the holomorphic structure without requiring the displayed smooth multipliers to be holomorphic in the initial trivialization.

Choose a positive Hermitian line bundle on \(S\) with integral first Chern form \(\eta\). Tensoring its pullback with \(\mathcal L\) gives the first Chern form \(2\beta+\pi^*\eta\). This form is positive definite: its two summands are positive respectively on the vertical and horizontal complex summands, and have no mixed terms. It is therefore an integral K\"ahler form on the compact manifold \(J\). Kodaira's embedding theorem implies that \(J\) is projective \cite[Main Theorem]{kodaira1954note}, \cite[\S3, Theorem~4]{kodaira1954}. Applying Chow's theorem to the graphs in products of projective varieties shows that the projection, the zero section, and the group operations are algebraic \cite[Theorem~VI, p.~912]{chow1949}. In particular \(\pi\) is an abelian scheme.

The zero section makes \(\pi_2(J)\to\pi_2(S)\) surjective and splits \(\pi_1(J)\to\pi_1(S)\). The homotopy exact sequence of the fibration therefore identifies the fibre lattice \(\Lambda^\vee\) with the full kernel on fundamental groups. Conjugation acts on this lattice by its given integral monodromy. Since \(h_*\) identifies \(\pi_1(S)\) with \(Q\), this proves the asserted semidirect product decomposition.
\end{proof}

\section{Proof of Theorem~A}\label{sec:whole-group}

We finish by embedding \(G\) as a subgroup of finite index in the fundamental group of an abelian scheme. Write \(A_{\mathbb Q}=A\otimes_{\mathbb Z}\mathbb Q\) and \(A_{\mathbb R}=A\otimes_{\mathbb Z}\mathbb R\).

\begin{lemma}\label{lem:integral}
Let \(1\to A\to G\xrightarrow{\rho}E\to1\) be an extension with \(G\) finitely generated and \(A\) free abelian of finite rank. Give \(A\) the \(E\)-action induced by conjugation. Suppose a crossed homomorphism \(t:G\to A_{\mathbb R}\) restricts on \(A\) to the natural inclusion. For some positive integer \(m\) there is an embedding \(\Phi:G\hookrightarrow A\rtimes E\) over \(E\) such that
\[
 \Phi(a)=(ma,1)\quad(a\in A),\qquad [A\rtimes E:\Phi(G)]=m^{\operatorname{rk}A}.
\]
\end{lemma}

\begin{proof}
Fix a finite generating set for \(G\) and a basis of \(A\). The values of a crossed homomorphism on the generators determine its value on any word by \(t(xy)=t(x)+\rho(x)t(y)\). The relations of \(G\), together with the prescribed values on words representing the basis of \(A\), impose rational affine linear equations on these generator values. Even when there are infinitely many relations, the augmented rows span a finite-dimensional rational vector space. A finite subsystem therefore has the same solutions. Since the system has a real solution, it has a rational solution \(t_{\mathbb Q}\). Choose \(m>0\) clearing the denominators on the generators. The action preserves \(A\), so \(mt_{\mathbb Q}(G)\subset A\). Now
\[
                    \Phi(g)=(m t_{\mathbb Q}(g),\rho(g))
\]
is a homomorphism. An element in its kernel belongs to \(A\) and is killed by \(m\), hence is trivial. The image surjects onto \(E\) and intersects \(A\) in \(mA\), which gives the stated index.
\end{proof}

\begin{proof}[Proof of Theorem~A]
Put \(G=\pi_1(M)\) and \(K=\operatorname{im}(\pi_1(F)\to G)\). Both are finitely generated, since \(M\) and \(F\) are compact manifolds. When \(B\) is a point, \(K=G\), and the theorem follows from the projective realization of virtually abelian Kähler groups \cite[Theorem~1.4]{baues2011}.

Suppose \(\dim B>0\). Apply Proposition~\ref{prop:neat} and write \(p:M^\circ\to\mathcal B^\circ\) for the resulting map. Its source still has fundamental group \(G\). Take a characteristic finite-index free abelian subgroup \(A\subset K\). Proposition~\ref{prop:transfer} gives a proper representable lift
\[
 p':M^\circ\longrightarrow\mathcal E^\circ,\qquad
 1\longrightarrow A\longrightarrow G\xrightarrow{\rho}
 E=\pi_1(\mathcal E^\circ)=\pi_1(\mathcal E)\longrightarrow1 .
\]
We use the analytification of the algebraic stack whenever the map is holomorphic. If \(A=0\), then \(G=E\) and Proposition~\ref{prop:projective-model} gives a smooth projective variety with fundamental group \(G\).

Suppose \(A\ne0\). For the open set \(U\) of Proposition~\ref{prop:neat}, put \(\mathcal V=\mathcal E|_U\), and let \(\widehat F\) be a geometric fibre of \(p'\) over \(\mathcal V\). It is compact, smooth and connected, and its fundamental group has image \(A\) in \(G\). Removing a proper analytic subset from \(M^\circ\) induces a surjection on fundamental groups. Applied to \((p')^{-1}(\mathcal V)\), this shows that \(\pi_1(\mathcal V)\to E\) is surjective. Let \(L\) denote its kernel. Proposition~\ref{prop:cocycle} gives an isomorphism of \(E\)-modules and a crossed homomorphism
\[
 H^1(\widehat F,\mathbb R)^L\simeq A_{\mathbb R}^{\vee},\qquad
 t:G\longrightarrow A_{\mathbb R},\qquad t(a)=a\quad(a\in A).
\]
Use \(A^\vee\) as the lattice in the invariant subspace. Its monodromy factors through \(E\), so this lattice defines a local system on all of \(\mathcal E\).

Over \(\mathcal V\), the local system \(R^1p'_*\mathbb R\) underlies a weight-one Hodge variation. If \(d=\dim\widehat F>0\), cup product with the restriction of a global Kähler class \([\omega]^{d-1}\) on the source, followed by fibre integration, gives the flat real polarization, with the weight-one sign convention. If \(d=0\), then \(A=0\), a case already settled. The open stack \(\mathcal V\) is a finite gerbe over the smooth quasi-projective variety \(U\), as required in Section~\ref{sec:hodge}.

By Proposition~\ref{prop:projective-model}, choose a morphism \(S\to\mathcal E\) from a smooth connected projective variety that induces an isomorphism on fundamental groups. Proposition~\ref{prop:jacobian} extends the invariant variation and constructs over \(S\) a projective abelian scheme \(J\to S\) with a zero section. Its first integral homology local system is determined by the lattice \((A^\vee)^\vee=A\) and the given \(E\)-action. Thus
\[
                         \pi_1(J)\simeq A\rtimes E .
\]
Lemma~\ref{lem:integral} realizes \(G\) as a subgroup of finite index in this group. Take the connected finite cover \(J'\to J\) corresponding to \(\Phi(G)\), and give it the pulled-back complex structure. The pullback of an integral Kähler class is an integral Kähler class, so \(J'\) is smooth projective. Since \(\pi_1(J')\simeq\Phi(G)\simeq G\), the variety \(J'\) realizes \(G\).
\end{proof}

\printbibliography

@incollection{arapura2011,
  shorthand   = {Ara11},
  author      = {Arapura, D.},
  title       = {Homomorphisms between {Kähler} groups},
  booktitle   = {Topology of Algebraic Varieties and Singularities},
  series      = {Contemporary Mathematics},
  number      = {538},
  publisher   = {American Mathematical Society},
  date        = {2011},
  pages       = {95--111},
  doi         = {10.1090/conm/538/10596},
}

@incollection{artinsga,
  shorthand   = {Art73},
  author      = {Artin, M.},
  title       = {Comparaison avec la cohomologie classique: cas d’un schéma lisse},
  booktitle   = {{SGA} 4, tome 3, exposé {XI}},
  series      = {Lecture Notes in Mathematics},
  number      = {305},
  publisher   = {Springer},
  date        = {1973},
  pages       = {64--78},
  doi         = {10.1007/BFb0070717},
}

@article{baues2011,
  shorthand    = {BR11},
  author       = {Baues, O. and Riesterer, J.},
  title        = {Virtually abelian {Kähler} and projective groups},
  journaltitle = {Abhandlungen aus dem Mathematischen Seminar der Universität Hamburg},
  volume       = {81},
  date         = {2011},
  pages        = {191--213},
  doi          = {10.1007/s12188-011-0056-1},
}

@article{cadman2007,
  shorthand    = {Cad07},
  author       = {Cadman, C.},
  title        = {Using stacks to impose tangency conditions on curves},
  journaltitle = {American Journal of Mathematics},
  volume       = {129},
  date         = {2007},
  pages        = {405--427},
  doi          = {10.1353/ajm.2007.0007},
}

@article{campana2004,
  shorthand    = {Cam04},
  author       = {Campana, F.},
  title        = {Orbifolds, special varieties and classification theory},
  journaltitle = {Annales de l’Institut Fourier},
  volume       = {54},
  number       = {3},
  date         = {2004},
  pages        = {499--630},
  doi          = {10.5802/aif.2027},
}

@article{camporb2011,
  shorthand    = {Cam11},
  author       = {Campana, F.},
  title        = {Orbifoldes géométriques spéciales et classification biméromorphe des variétés kählériennes compactes},
  journaltitle = {Journal of the Institute of Mathematics of Jussieu},
  volume       = {10},
  date         = {2011},
  pages        = {809--934},
  doi          = {10.1017/S1474748010000101},
  eprinttype   = {arxiv},
  eprint       = {0705.0737v8},
}

@article{cce2014,
  shorthand    = {CCE14},
  author       = {Campana, F. and Claudon, B. and Eyssidieux, P.},
  title        = {Représentations linéaires des groupes kählériens et de leurs analogues projectifs},
  journaltitle = {Journal de l’École polytechnique — Mathématiques},
  volume       = {1},
  date         = {2014},
  pages        = {331--342},
  doi          = {10.5802/jep.12},
}

@unpublished{cceadd2015,
  shorthand = {CCE15},
  author    = {Campana, F. and Claudon, B. and Eyssidieux, P.},
  title     = {Addendum à l’article « Représentations linéaires des groupes kählériens et de leurs analogues projectifs »},
  date      = {2015-05-22},
  url       = {https://claudon.pages.math.cnrs.fr/pdf/addendum_JEP.pdf},
}

@article{claudon2019,
  shorthand    = {CHL19},
  author       = {Claudon, B. and Höring, A. and Lin, H.-Y.},
  title        = {The fundamental group of compact {Kähler} threefolds},
  journaltitle = {Geometry \& Topology},
  volume       = {23},
  number       = {7},
  date         = {2019},
  pages        = {3233--3271},
  doi          = {10.2140/gt.2019.23.3233},
}

@article{chow1949,
  shorthand    = {Cho49},
  author       = {Chow, W.-L.},
  title        = {On compact complex analytic varieties},
  journaltitle = {American Journal of Mathematics},
  volume       = {71},
  number       = {4},
  date         = {1949},
  pages        = {893--914},
  doi          = {10.2307/2372375},
}

@article{claudon2018,
  shorthand    = {Cla18},
  author       = {Claudon, B.},
  title        = {Smooth families of tori and linear {Kähler} groups},
  journaltitle = {Annales de la Faculté des Sciences de Toulouse. Mathématiques},
  series       = {6},
  volume       = {27},
  number       = {3},
  date         = {2018},
  pages        = {477--496},
  doi          = {10.5802/afst.1576},
}

@article{deligne1968,
  shorthand    = {Del68},
  author       = {Deligne, P.},
  title        = {Théorème de {Lefschetz} et critères de dégénérescence de suites spectrales},
  journaltitle = {Publications Mathématiques de l'IHÉS},
  volume       = {35},
  date         = {1968},
  pages        = {107--126},
  doi          = {10.1007/BF02698925},
}

@incollection{deligne1987,
  shorthand = {Del87},
  author    = {Deligne, P.},
  title     = {Un théorème de finitude pour la monodromie},
  booktitle = {Discrete Groups in Geometry and Analysis},
  series    = {Progress in Mathematics},
  number    = {67},
  publisher = {Birkhäuser},
  date      = {1987},
  pages     = {1--19},
  url       = {https://publications.ias.edu/sites/default/files/56_Untheoremede.pdf},
}

@incollection{deligne1994,
  shorthand = {Del94},
  author    = {Deligne, P.},
  title     = {Décompositions dans la catégorie dérivée},
  booktitle = {Motives},
  series    = {Proceedings of Symposia in Pure Mathematics},
  number    = {55, part 1},
  date      = {1994},
  pages     = {115--128},
  url       = {https://publications.ias.edu/sites/default/files/68_Decompositionsdans.pdf},
}

@unpublished{dejong,
  shorthand = {dJ},
  author    = {de Jong, A. J.},
  title     = {A result of {Gabber}},
  note      = {Author manuscript},
  url       = {https://www.math.columbia.edu/~dejong/papers/2-gabber.pdf},
}

@article{ehkv2001,
  shorthand    = {EHKV01},
  author       = {Edidin, D. and Hassett, B. and Kresch, A. and Vistoli, A.},
  title        = {{Brauer} groups and quotient stacks},
  journaltitle = {American Journal of Mathematics},
  volume       = {123},
  date         = {2001},
  pages        = {761--777},
  doi          = {10.1353/ajm.2001.0024},
}

@thesis{giraud1966,
  shorthand   = {Gir66},
  author      = {Giraud, J.},
  title       = {Cohomologie non abélienne de degré 2},
  type        = {phdthesis},
  institution = {Faculté des sciences de l'Université de Paris},
  date        = {1966},
  url         = {https://archive.org/details/IA41544320_0022},
}

@unpublished{huwu2026,
  shorthand  = {HW26},
  author     = {Hu, W. and Wu, P.},
  title      = {Projective surfaces with infinite fundamental group and trivial profinite completion},
  date       = {2026},
  note       = {Preprint},
  eprinttype = {arxiv},
  eprint     = {2611.17847v1},
}

@incollection{johnson1991,
  shorthand = {Joh91},
  author    = {Johnson, F. E. A.},
  title     = {Flat algebraic manifolds},
  booktitle = {Geometry of Low-Dimensional Manifolds, 1},
  series    = {London Mathematical Society Lecture Note Series},
  number    = {150},
  publisher = {Cambridge University Press},
  date      = {1991},
  pages     = {73--92},
  doi       = {10.1017/CBO9780511629334.007},
}

@article{johnsonrees1987,
  shorthand    = {JR87},
  author       = {Johnson, F. E. A. and Rees, E. G.},
  title        = {On the fundamental group of a complex algebraic manifold},
  journaltitle = {Bulletin of the London Mathematical Society},
  volume       = {19},
  date         = {1987},
  pages        = {463--466},
  doi          = {10.1112/blms/19.5.463},
}

@article{kodaira1954note,
  shorthand    = {Kod54a},
  author       = {Kodaira, K.},
  title        = {On {Kähler} varieties of restricted type},
  journaltitle = {Proceedings of the National Academy of Sciences of the United States of America},
  volume       = {40},
  date         = {1954},
  pages        = {313--316},
  doi          = {10.1073/pnas.40.5.313},
}

@article{kodaira1954,
  shorthand    = {Kod54b},
  author       = {Kodaira, K.},
  title        = {On {Kähler} varieties of restricted type (an intrinsic characterization of algebraic varieties)},
  journaltitle = {Annals of Mathematics},
  series       = {2},
  volume       = {60},
  number       = {1},
  date         = {1954},
  pages        = {28--48},
  doi          = {10.2307/1969701},
}

@article{kv2004,
  shorthand    = {KV04},
  author       = {Kresch, A. and Vistoli, A.},
  title        = {On coverings of {Deligne--Mumford} stacks and surjectivity of the {Brauer} map},
  journaltitle = {Bulletin of the London Mathematical Society},
  volume       = {36},
  date         = {2004},
  pages        = {188--192},
  doi          = {10.1112/S0024609303002728},
}

@article{lin2020,
  shorthand    = {Lin20},
  author       = {Lin, H.-Y.},
  title        = {Algebraic approximations of compact {Kähler} manifolds of algebraic codimension 1},
  journaltitle = {Duke Mathematical Journal},
  volume       = {169},
  number       = {14},
  date         = {2020},
  pages        = {2767--2826},
  doi          = {10.1215/00127094-2020-0018},
}

@unpublished{linthreefolds,
  shorthand  = {Lin24},
  author     = {Lin, H.-Y.},
  title      = {Algebraic approximations of compact {Kähler} threefolds},
  date       = {2024},
  eprinttype = {arxiv},
  eprint     = {1710.01083v4},
}

@unpublished{noohi2005,
  shorthand  = {Noo05},
  author     = {Noohi, B.},
  title      = {Foundations of {Topological Stacks I}},
  note       = {Preprint},
  date       = {2005},
  eprinttype = {arxiv},
  eprint     = {math/0503247v1},
}

@unpublished{openai2026,
  shorthand = {OAI26a},
  author    = {{OpenAI}},
  title     = {The abelianity conjecture for special compact {Kähler} manifolds},
  note      = {Preprint},
  date      = {2026},
  url       = {https://github.com/openai/math/blob/3344da59991ed94f9ba711217b29c1d760badfa6/preprints/The-abelianity-conjecture-for-special-compact-Kahler-manifolds-September-23-2026/paper.pdf},
}

@unpublished{openai2026b,
  shorthand = {OAI26b},
  author    = {{OpenAI}},
  title     = {A projective fourfold with large fundamental group and non-{Stein} universal cover},
  note      = {Preprint},
  date      = {2026},
  url       = {https://github.com/openai/math/blob/3344da59991ed94f9ba711217b29c1d760badfa6/preprints/A-projective-fourfold-with-large-fundamental-group-and-non-Stein-universal-cover-October-5-2026/large-fundamental-group-non-stein.pdf},
}

@article{saito1990,
  shorthand    = {Sai90},
  author       = {Saito, M.},
  title        = {Decomposition theorem for proper {Kähler} morphisms},
  journaltitle = {Tohoku Mathematical Journal},
  volume       = {42},
  date         = {1990},
  pages        = {127--148},
  doi          = {10.2748/tmj/1178227650},
}

@unpublished{saito2022,
  shorthand  = {Sai22},
  author     = {Saito, M.},
  title      = {Some remarks on decomposition theorem for proper {Kähler} morphisms},
  note       = {Preprint},
  date       = {2022},
  eprinttype = {arxiv},
  eprint     = {2204.09026v5},
}

@article{schmid1973,
  shorthand    = {Sch73},
  author       = {Schmid, W.},
  title        = {Variation of {Hodge} structure: the singularities of the period mapping},
  journaltitle = {Inventiones Mathematicae},
  volume       = {22},
  date         = {1973},
  pages        = {211--319},
  doi          = {10.1007/BF01389674},
}

\bigskip
\noindent\begin{minipage}{\textwidth}
\small
\textsc{Tianshu Jiang}\\
Nanjing University\\
\href{mailto:jts2026@nju.edu.cn}{\nolinkurl{jts2026@nju.edu.cn}}

\medskip
\textsc{Jiayu Li}\\
Nanjing University\\
\href{mailto:jiayuli@nju.edu.cn}{\nolinkurl{jiayuli@nju.edu.cn}}
\end{minipage}
\end{document}